\documentclass[12pt,notitlepage, reqno]{amsart}
\usepackage{amsmath}
\usepackage{amsfonts,amsthm,amssymb}
\usepackage{euscript}
\makeatletter
 \def\activeat#1{\csname @#1\endcsname}
 \def\def@#1{\expandafter\def\csname @#1\endcsname}
 {\catcode`\@=\active \gdef@{\activeat}}

\let\ssize\scriptstyle
\newdimen\ex@	\ex@.2326ex

 \def\requalfill{\cleaders\hbox{$\mkern-2mu\mathord=\mkern-2mu$}\hfill
  \mkern-6mu\mathord=$}
 \def\eqfill{$\m@th\mathord=\mkern-6mu\requalfill}
 \def\deffill{\hbox{$:=$}$\m@th\mkern-6mu\requalfill}
 \def\fiberbox{\hbox{$\vcenter{\hrule\hbox{\vrule\kern1ex
     \vbox{\kern1.2ex}\vrule}\hrule}$}}

 \newdimen\arrwd 
 \newdimen\minCDarrwd \minCDarrwd=2.5pc
 	
 \def\findarrwd#1#2#3{\arrwd=#3%
  \setbox\z@\hbox{$\ssize\;{#1}\;\;$}%
  \setbox\@ne\hbox{$\ssize\;{#2}\;\;$}%
  \ifdim\wd\z@>\arrwd \arrwd=\wd\z@\fi
  \ifdim\wd\@ne>\arrwd \arrwd=\wd\@ne\fi}
 \newdimen\arrowsp\arrowsp=0.375em  	
 \def\findCDarrwd#1#2{\findarrwd{#1}{#2}{\minCDarrwd}
    \advance\arrwd by 2\arrowsp}
 \newdimen\minarrwd 
 \setbox\z@\hbox{$\longrightarrow$} \minarrwd=\wd\z@

 \def\harrow#1#2#3#4{{\minarrwd=#1\minarrwd
   \findarrwd{#2}{#3}{\minarrwd}\kern\arrowsp
    \mathrel{\mathop{\hbox to\arrwd{#4}}\limits^{#2}_{#3}}\kern\arrowsp}}

 \def@]#1>#2>#3>{\harrow{#1}{#2}{#3}\rightarrowfill}
 \def@>#1>#2>{\harrow1{#1}{#2}\rightarrowfill}
 \def@<#1<#2<{\harrow1{#1}{#2}\leftarrowfill}
 \def@={\harrow1{}{}\eqfill}
 \def@:#1={\harrow1{}{}\deffill}
 \def@ N#1N#2N{\vCDarrow{#1}{#2}\UpDownarrow}
 \def\UpDownarrow{\uparrow\,\Big\downarrow}

\def\hookrightarrowfill{\hbox{$\lhook\joinrel$}\rightarrowfill}
\def@{hk>}#1>#2>{\harrow1{#1}{#2}\hookrightarrowfill}
\def@ c>#1>#2>{\harrow1{#1}{#2}\hookrightarrowfill}
\def\hookleftarrowfill{\leftarrowfill\hbox{$\joinrel\rhook$}}
\def@{hk<}#1<#2<{\harrow1{#1}{#2}\hookleftarrowfill}

 \def@.{\ifodd\row\relax\harrow1{}{}\hfill
   \else\vCDarrow{}{}.\fi}
 \def@|{\vCDarrow{}{}\Vert}
 \def@ V#1V#2V{\vCDarrow{#1}{#2}\downarrow}
 \def@ A#1A#2A{\vCDarrow{#1}{#2}\uparrow}
 \def@(#1){\arrwd=\csname col\the\col\endcsname\relax
   \hbox to 0pt{\hbox to \arrwd{\hss$\vcenter{\hbox{$#1$}}$\hss}\hss}}

 \def\squash#1{\setbox\z@=\hbox{$#1$}\finsm@@sh}
\def\finsm@@sh{\ifnum\row>1\ht\z@\z@\fi \dp\z@\z@ \box\z@}

 \newcount\row \newcount\col \newcount\numcol \newcount\arrspan
 \newdimen\vrtxhalfwd  \newbox\tempbox

 \def\innernewdimen{\alloc@1\dimen\dimendef\insc@unt}
 \def\measureinit{\col=1\vrtxhalfwd=0pt\arrspan=1\arrwd=0pt 
   \setbox\tempbox=\hbox\bgroup$}
 \def\setinit{\col=1\hbox\bgroup$\ifodd\row
   \kern\csname col1\endcsname
   \kern-\csname row\the\row col1\endcsname\fi}
 \def\findvrtxhalfsum{$\egroup
  \expandafter\innernewdimen\csname row\the\row col\the\col\endcsname
  \global\csname row\the\row col\the\col\endcsname=\vrtxhalfwd
  \vrtxhalfwd=0.5\wd\tempbox
  \global\advance\csname row\the\row col\the\col\endcsname by \vrtxhalfwd 
  \advance\arrwd by \csname row\the\row col\the\col\endcsname
  \divide\arrwd by \arrspan
  \loop\ifnum\col>\numcol \numcol=\col%
     \expandafter\innernewdimen \csname col\the\col\endcsname
     \global\csname col\the\col\endcsname=\arrwd
   \else \ifdim\arrwd >\csname col\the\col\endcsname
      \global\csname col\the\col\endcsname=\arrwd\fi\fi
   \advance\arrspan by -1 %
   \ifnum\arrspan>0 \repeat}
 \def\setCDarrow#1#2#3#4{\advance\col by 1 \arrspan=#1 
    \arrwd= -\csname row\the\row col\the\col\endcsname\relax
    \loop\advance\arrwd by \csname col\the\col\endcsname
     \ifnum\arrspan>1 \advance\col by 1 \advance\arrspan by -1%
     \repeat
    \squash{\mathop{
     \hbox to\arrwd{\kern\arrowsp#4\kern\arrowsp}}\limits^{#2}_{#3}}}
 \def\measureCDarrow#1#2#3#4{\findvrtxhalfsum\advance\col by 1%
   \arrspan=#1\findCDarrwd{#2}{#3}%
    \setbox\tempbox=\hbox\bgroup$}
 \def\vCDarrow#1#2#3{\kern\csname col\the\col\endcsname
    \hbox to 0pt{\hss$\vcenter{\llap{$\ssize#1$}}%
     \Big#3\vcenter{\rlap{$\ssize#2$}}$\hss}\advance\col by 1}

 \def\setCD{\def\harrow{\setCDarrow}%
  \def\\{$\egroup\advance\row by 1\setinit}
  \m@th\lineskip3\ex@\lineskiplimit3\ex@ \row=1\setinit}
 \def\endsetCD{$\egroup}
 \def\drop#1\\{\findvrtxhalfsum\advance\row by 2 \measureinit}
 \def\measure{\bgroup
  \def\harrow{\measureCDarrow}%
  \def\\##1{\ifx##1\endmeasure\endmeasure\else\expandafter\drop\fi}%
  \row=1\numcol=0\measureinit}
 \def\endmeasure{\findvrtxhalfsum\egroup}

 \newcount\savedcount	
 \def\LCD#1\end{\savedcount=\count11
   \measure#1\endmeasure
   \vcenter{\setCD#1\endsetCD\kern\medskipamount}%
   \global\count11=\savedcount\end}

 \newenvironment{CD}{\let\at=@\catcode`\@=\active\LCD}{\catcode`\@=12\relax}
\makeatother

\advance\voffset by -0.4 truein
\advance\hoffset by -0.6 truein
\font\tenbi=cmmi14
 \font\sevenbi=cmmi10 \font\fivebi=cmmi7
 \newfam\bifam  \textfont\bifam=\tenbi
 \scriptfont\bifam=\sevenbi \scriptscriptfont\bifam=\fivebi
 \mathchardef\variablemega="7121  
 \mathchardef\variablenu="7117 

\font\smallrm=cmr8

\renewcommand{\c}{\mathcal}
\newcommand{\Og}{\Omega}
\newcommand{\T}{{\EuScript{T}}}

\newcommand{\h}{\operatorname{h}}
\newcommand{\Hom}{\operatorname{Hom}}
\newcommand{\Ext}{\operatorname{Ext}}
\newcommand{\Cok}{\operatorname{Cok}}
\newcommand{\Ker}{\operatorname{Ker}}
\newcommand{\IM}{\operatorname{Im}}
\newcommand{\rk}{\operatorname{rk}}
\newcommand{\reg}{\operatorname{reg}}
\newcommand{\Sing}{\operatorname{Sing}}
\newcommand{\tr}{\operatorname{tr}}
\newcommand{\wt}{\overline}
\renewcommand{\:}{\colon}

\newcommand{\ox}{\otimes}
\newcommand{\IP}{\text{\bf P}}
\def\I#1{\mathbf#1}

\newcommand{\bigomega}{\hbox{\large $\omega$}}
\newcommand{\w}{\bigomega}
\def\fld #1#2#3{#1\,\frac{\partial}{\partial x}
	+#2\,\frac{\partial}{\partial y}
	+#3\,\frac{\partial}{\partial z}}
\def\Ham #1#2#3{\frac{\partial#1}{\partial#2}\,\frac{\partial}{\partial#3}
 	-\frac{\partial#1}{\partial#3}\,\frac{\partial}{\partial#2}}
\let~=\thinspace

\newtheorem{theorem}{Theorem}[section]
\newtheorem{corollary}[theorem]{Corollary}
\newtheorem{lemma}[theorem]{Lemma}
\newtheorem{proposition}[theorem]{Proposition}
\theoremstyle{definition}
\newtheorem{subsct}[theorem]{}
\theoremstyle{plain}
 
\begin{document}

\title[{\smallrm Bounds on leaves of one-dimensional foliations}]
      {Bounds on leaves of one-dimensional foliations}
\author[{\smallrm E. Esteves and S. Kleiman}]{E. Esteves$^1$
\ and \ S. Kleiman}
 \date{10 September 2002}

\begin{abstract} 
Let $X$ be a variety over an algebraically closed field,
$\eta\:\Omega^1_X\to\c L$ a one-dimensional singular foliation, and
$C\subseteq X$ a projective leaf of $\eta$.  We prove that
$2p_a(C)-2=\deg(\c L|C)+\lambda(C)-\deg(C\cap S)$ where $p_a(C)$ is the
arithmetic genus, where $\lambda(C)$ is the colength in the dualizing
sheaf of the subsheaf generated by the K\"ahler differentials, and where
$S$ is the singular locus of $\eta$.  We bound $\lambda(C)$ and
$\deg(C\cap S)$, and then improve and extend some recent results of
Campillo, Carnicer, and de la Fuente, and of du Plessis and Wall.
 \end{abstract}
\maketitle

\section{Introduction}

In 1891, Poincar\'e \cite{P}, p.~161, considered, in effect, a foliation
of the plane given by a polynomial vector field, and he posed the
problem of deciding whether it is algebraically integrable or not.
Poincar\'e observed that it is enough to find a bound on the degree of
the integral.

Over the years, this problem has attracted a lot of attention.  
Recently, it has been interpreted as the problem of bounding the degrees of 
the algebraic leaves of the foliation, be it algebraically integrable or 
not. As such, the problem was
addressed in \cite{CLN}, by local methods, and in \cite{CC}, \cite{CO}, 
and \cite{Lc}, using resolution of singularities. A bound depending only 
on the degree of the foliation was proved in
\cite{C} for foliations without diacritical singularities.

In general, Lins Neto \cite{LN}, Main Thm., p.~234,  showed that there is
no bound depending only on the degree of the foliation and on the analytic
type of its singularities. Bounds depending on the degree of the 
foliation and the analytic type of the singularities of the leaves were 
proved in \cite{CL}, \cite{PW} and \cite{W}. In \cite{Pe}, bounds depending 
on the degree and plurigenera of the foliation and the geometric genera of 
the leaves were proved for foliations of general type.

The problem was extended to surfaces with trivial Picard group in 
\cite{B} and, more generally, to any smooth ambient variety in \cite{CCF}. 
Bounds on numerical invariants of subvarieties saturated by leaves were 
considered in \cite{E}, \cite{So1} and \cite{So2}. Finally, the analogous 
problem for Pfaff differential equations was considered in \cite{BM} and 
\cite{EK}.
\renewcommand{\thefootnote}{}%
\footnote
 {2000 {\it Mathematics Subject Classification} 37F75 (primary), 14H50,
   32S65, 14H20 (secondary).}
\footnote
 {{\it Keywords} Foliations, curves, singularities.}
\footnote
 {$^1$This author is grateful to A. Campillo, L. G. Mendes, P. Sad,
M. Soares, and especially J. V. Pereira for helpful discussions on the
subject matter.  He is also grateful to CNPq for a grant, 
Proc.\ 202151/90-5, supporting
a year-long visit to MIT, and grateful to MIT for its hospitality. 
He was also supported by PRONEX, Conv\^enio 41/96/0883/00, CNPq, 
Proc.\ 300004/95-8, and FAPERJ, Proc.{} E-26/170.418/2000-APQ1.}

Here we address the following version of the problem.  Let $X$ be a
variety over an algebraically closed field of arbitrary
characteristic.  Let $C\subseteq X$ be a {\it curve}, that is, a {\bf
  reduced} subscheme of pure dimension 1; assume $C$ is projective.  Let
$\eta\:\Og^1_X\to\c L$ be a {\it  one-dimensional singular foliation}
of $X$; that is, $\eta$ is nonzero, and $\c L$ is invertible.  Assume
$C$ is a \emph{leaf}; that is, $C$ contains only finitely many
singularities of $\eta$, and $\eta|C$ factors through the standard map
$\sigma\:\Og^1_X|C\to\Og^1_C$.  Say $\mu\:\Og^1_C\to\c L|C$ is the
induced map.  We strive to relate the numerical invariants of $C$ and
$\mu$.

The major global  invariant of $C$ is its arithmetic genus,
$p_a(C):=1-\chi(\c O_C)$.  Notice that $p_a(C)=h^1(\c O_C)$ when  $C$ is
connected, and that $p_a(C)$ remains constant when $C$ varies in a family.

The singularities $P$ of $C$ are measured by several
invariants.  One in particular arises naturally in the present work.  It
is denoted $\lambda(C,P)$ by Buchweitz and Greuel in \cite{BG},
Def.~6.1.1, p.~265, and it is defined as the colength, in the dualizing
module $\w_P$, of the $\c O_{C,P}$-submodule generated by $\Og^1_{C,P}$.
Notice that $\lambda(C,P)>0$ if and only if $P$ is singular.  So we may
set $\lambda(C):=\sum\lambda(C,P)$.

Our key relation is the following simple formula, given in Proposition
5.2:
 $$
 2p_a(C)-2-\deg(\c L|C)=\lambda(C)-\deg(C\cap S)\eqno(1.1).
 $$
Here $S$ is the \emph{singular locus} of $\eta$, that is, the subscheme
of $X$ where $\eta$ fails to be surjective; so $C\cap S$ is the singular
locus of $\mu$.  We prove our formula by comparing Euler characteristics
of certain torsion-free sheaves on $C$.

Under more restrictive hypotheses, versions of Formula (1.1) were
proved by Cerveau and Lins Neto \cite{CLN}, Prop., p.~885, and by Lins
Neto and Soares \cite{LS}, Prop.~2.7, p.~659.  In \cite{So1},
p.~495, Soares suggested using the formula when $C$ is smooth, to
solve the Poincar\'e problem by bounding $\deg(C\cap S)$ from below.
In the same vein, our main results, Theorem 5.3 and Theorem 6.1,
follow from the general case of Formula (1.1) and from bounds we
obtain on $\lambda(C,P)$ and $\deg(C\cap S)$.

Note that $\deg(C\cap S)\geq \iota(C)$ where $\iota(C):=\sum\iota(C,P)$
and $\iota(C,P)$ is the least length of the cokernel of a map
$\Og^1_{C,P} \to \c O_{C,P}$.  Hence, as is also stated in our
Proposition 5.2,
 $$2p_a(C)-2-\deg(\c L|C)\le\lambda(C)-\iota(C).\eqno(1.2)$$
In characteristic 0, if $P$ is a singularity, then $\iota(C,P)\ge1$
because $\Og^1_{C,P}/\text{torsion}$ cannot be free by \cite{L},
Thm.~1, p.~879.  Hence then $\iota(C)$ is at least the number of
singularities.

Assume $X$ is smooth.  In \cite{CCF}, Thm.~2.7, p.~62, Campillo,
Carnicer and De la Fuente gave an upper bound on $2p_a(C)-2-\deg(\c
L|C)$ in terms of multiplicities associated to $C$ and $\eta$ along a
sequence of blowups of $X$ resolving the singularities of $C$. As a
consequence, they obtained in \cite{CCF}, Thm.~3.1, p.~64, an upper bound
on $2p_a(C)-2-\deg(\c L|C)$ that holds universally for all
$\eta$ having $C$ as leaf.  Our Theorem 5.3 provides a somewhat better
bound; this bound follows from (1.2), given the bound on $\lambda(C,P)$
asserted in our Proposition 4.4.  Thus (1.2) is the sharpest available
bound on $2p_a(C)-2-\deg(\c L|C)$.

Our proof of Proposition 4.4 uses the Hironaka--Noether bound,
Proposition 3.1.  It bounds the colength $\ell$ of a reduced
one-dimensional Noetherian local ring $A$ in the blowup at its maximal
ideal in terms of its multiplicity $e$; namely, $\ell\leq e(e-1)/2$,
with equality if and only if $A$ has embedding dimension at most 2.
Noether \cite{N} considered, in effect, only the case where $A$ is the
local ring of a complex plane curve.  Hironaka \cite{Hi}, p.~186,
asserted the bound without proof when $A$ is the local ring of an
arbitrary complex curve.  In the same setup, Stevens \cite{St} proved a
formula for $\ell$, and then asserted the bound without proof.  Inspired
by the Stevens's work, we give a somewhat different proof, and obtain
the general case.

Take $X:=\IP^n$ now, and set $d:=\deg C$.  Suppose $d$ is not a multiple
of the characteristic.  Over $\mathbb C$, Jouanolou \cite{J}, Prop.~4.2,
p.~130, proved $C\cap S$ is nonempty, even when $C$ is smooth. In
\cite{EK}, Cor.~4.5, Jouanolou's result is refined: the
Castelnuovo--Mumford regularity $\reg(C\cap S)$ is shown to be at least
$m+1$ where $m:=1+\deg\c L$.  Now, the regularity of any finite
subscheme is at least its degree.  Hence, (1.1) yields
	$$2p_a(C)-(d-1)(m-1)\le\lambda(C),\eqno(1.3)$$
 which our Theorem 6.1 asserts.  It continues by asserting that, if
equality holds, then $\deg(C\cap S)=m+1$; also then $C\cap S$ lies on a
line $M$, and either $M\subseteq S$ or $M$ is a leaf.

If, in addition, the singular locus $S$ is finite, then, as our
Proposition 6.3 asserts,
	$$\lambda(C)\le2p_a(C)-(d-1)(m-1)+m^2+\dotsb+m^n.\eqno(1.4)$$
This bound too results from (1.1); indeed, a simple Chern class
computation evaluates $\deg(S)$, but  $\deg(S)\ge\deg(C\cap S)$.

Another important global invariant of $C$ is its geometric genus,
$p_g(C):=\h^1(\c O_{\wt C})$ where $\wt C$ is the normalization of $C$.
Our Corollary 6.2 asserts that, if $C$ is connected and the
characteristic is $0$, then
 $$p_g(C)\leq (m-1)(d-1)/2+(r(C)-1)/2$$
 where $r(C)$ is the number of irreducible components.
Notice that this bound is nontrivial for $m<d-1$ and that it does not
depend in any way on the singularities of $C$ or of $\eta$.  The
problem of bounding $p_g(C)$ was posed by Painlev\'e and has been
considered by Lins Neto among others; see \cite{LN}.

There are two better known singularity invariants, the
$\delta$-invariant $\delta(C,P)$ and the Tjurina number $\tau(C,P)$.
The former is the colength of $\c O_{C,P}$ in its normalization; the
latter, the dimension of the tangent space of the miniversal deformation
space of the singularity.  These invariants are related to
$\lambda(C,P)$.  First,
 $\delta(C,P)\leq\lambda(C,P)\leq 2\delta(C,P),$
 but the second inequality is valid only in characteristic 0; see
 Subsection 2.1.  Second, $\tau(C,P)=\lambda(C,P)$ if $C$ is a complete
 intersection at $P$; see Proposition 2.2.

Finally, take $X:=\IP^2$.  Then $p_a(C)=(d-1)(d-2)/2$, and
$\lambda(C)=\tau(C)$ where $\tau(C):=\sum\tau(C,P)$.  Again suppose $d$
is not a multiple of the characteristic.  Then (1.3) and (1.4) hold, and
reduce to the following lower and upper bounds on $\tau(C)$:
	$$(d-1)(d-m-1)\leq\tau(C)\leq (d-1)(d-m-1)+m^2.\eqno(1.5)$$

These bounds are 
the ones masterfully proved over $\mathbb C$
by du Plessis and Wall \cite{PW}, Thm.~3.2, p.~263, in a more
elementary way.  However, they define $m$ as the least degree of a
nontrivial polynomial vector field $\phi$ annihilating the equation of
$C$.  
Considering the foliation $\eta$ defined by $\phi$, we derive their 
lower bound in our Corollary 6.4. Their upper bound is also obtained 
there, under the additional assumption that the singular locus of 
$\eta$ intersects $C$ in finitely many points.

In fact, Du Plessis and Wall prove more: if $2m+1> d$, then
	$$\tau(C)\leq (d-1)(d-m-1)+m^2-(2m+2-d)(2m+1-d)/2.$$
 The present authors have found a more conceptual version of the
proof, which also works for other ambient spaces; this material is
treated in \cite{EK2}.

The lower bound in (1.5) was rediscovered over $\mathbb C$ by Chavarriga
and Llibre \cite{CL}, Thm.~3, p.~12, and they gave yet a third proof.

The lower bound in (1.5) is improved in characteristic 0 via yet a
fourth argument in \cite{EK2}, as follows:
$(d-1)(d-m-1)+u\leq\tau(C)$ where $u$ is the number of singularities
{\it not\/} quasi-homogeneous (that is, at which a local analytic
equation is not weighted homogeneous); moreover, if equality holds, then
either $m=d-1$ and $C$ is smooth, or $m<d-1$ and $\reg(\Sing C)=2d-3-m$.

The Poincar\'e problem is to bound $d$ given the invariants of
$\eta$.  As is well known, the difficulty lies in the possibility that
$C$ may be highly singular.  In this connection, the lower bound in
(1.5) says this: the higher its degree, the more singular is $C$.  As
noted above, our proof of (1.5) uses the lower bound $\reg(C\cap S)$
given in \cite{EK}, Cor.~4.5.  A result in \cite{EK2} 
asserts that $\reg(\Sing C)\ge2d-3-m$ if $m \le d-2$ and that $\reg(\Sing
C)=2d-3-m$ if $m\le(d-2)/2$, provided $d$ is not a multiple of the
characteristic.  In other words, for high $d$, not only must $C$ have
many singularities, but also they must lie in special position in the
plane.

In short, Section 2 of the present paper introduces some local and some
global invariants of a curve $C$, and relates them.  Section 3 treats
the Hironaka--Noether bound.  Section 4 uses this bound to help
establish an upper bound on $\lambda(C,P)$.  Section 5 establishes our
bound (1.2) on $2p_a(C)-2-\deg(\c L|C)$, and compares it favorably to
the bound of Campillo, Carnicer and De la Fuente with the aid of our
bound on $\lambda(C,P)$.  Finally, Section 6 establishes the bounds
(1.3) and (1.4) on $\lambda(C)$, and shows that they recover the bounds
in (1.5) on $\tau(C)$ in the form treated by du Plessis and Wall and by
Chavarriga and Llibre.

\section{Invariants of  curves}

\begin{subsct} \emph{Local invariants.} 
Let $C$ be a curve, $\I n\:\wt C\to C$ the normalization map, and
	$$\I n^\#\:\c O_C\to\I n_*\c O_{\wt C}\text{ \ and \ }
	d\,\I n\:\Og^1_C\to\I n_*\Og^1_{\wt C}$$
 the associated maps on sheaves of functions and differentials.  Let
$\w_C$
 be the dualizing sheaf (or canonical sheaf, or Rosenlicht's sheaf of
 regular differentials); see  \cite{Ro}, or \cite{Ha}, Sec.~III-7, or
 \cite{BG}, pp.~243--244, or \cite{AK}, for example.  There is a natural
map
	$$\tr\:\I n_*\Og^1_{\wt C}\to\w_C;$$
 it is known as the {\it trace}, and the composition
	$$\begin{CD}
   \gamma\:\Og^1_C @>d\,\I n>>\I n_*\Og^1_{\wt C}@>\tr>>\w_C
	\end{CD}$$
 is  known as the {\it class map}.

Fix a closed point $P\in C$.  Taking lengths
$\ell(-)$, set \begin{align*}
 &\delta(C,P):=\ell(\Cok(\I n^\#_P)),\\
 &\tau(C,P):=\ell(\Ext^1_{\c O_{C,P}}(\Og^1_{C,P},\c O_{C,P})),\\
 &\lambda(C,P):=\ell(\Cok(\gamma_P)).
\end{align*}
 The first two invariants are known respectively as the {\it
  $\delta$-invariant} or the {\it genus  diminution,} and the
  {\it Tjurina number;} see \cite{St}, p.~98, and \cite{G},
  pp.~142--143.  The third invariant was formally introduced and studied
  by Buchweitz and Greuel \cite{BG}, pp.~265--269, although it appears
  implicitly earlier, notably in Rim's paper \cite{R}. 

By Rosenlicht's theorem (see \cite{Ro}, Thm.~8 and Cor.~1, pp.~177--178,
or \cite{AK}, Prop.~1.16(ii), p.~168), the cokernels of $\I n^\#$ and
$\tr$ are perfectly paired; so $\delta(C,P)=\ell(\Cok(\tr_P))$.  Hence
 $$
 \lambda(C,P)=\delta(C,P)+\ell(\Cok(d\,\I n_P)).\eqno(2.1.1)
 $$

Let $\alpha\:\Og^1_{C,P}\to\c O_{C,P}$ range over all maps such that
$\Cok\alpha$ has finite length, and set
	$$\iota(C,P):=\min_\alpha\ell(\Cok\alpha).$$
This invariant is the local isomorphism defect of
$\Og^1_C/\text{torsion}$ in $\c O_C$ at $P$, as defined by Greuel and
Lossen in \cite{GL}, p.~330, and as defined earlier, but with the
opposite sign, by Greuel and Karras in \cite{GrK}, p.~103; however,
the present invariant $\iota(C,P)$ itself is not explicitly considered
in either of those papers.

Suppose that $P$ is a singularity of $C$.  In characteristic zero,
$\iota(C,P)\ge1$ because $\Hom(\Og^1_{C,P},\c O_{C,P})$ is not free by
\cite{L}, Thm.~1, p.~879.  In characteristic $p>0$, sometimes
$\iota(C,P)=0$; for example (see \cite{L}, p.~892), in the plane, take
$C:y^{p+1}-x^p=0$ and take $P:=(0,0)$.

Let $r(C,P)$ denote the number of branches, or analytic components, of
$C$ at $P$.

Let $d\:\c O_C\to\Og^1_C$ be the universal derivation, and set
 \[
 \mu(C,P):=\ell(\Cok(\gamma\circ d)_P).
 \]
Then $\lambda(C,P)\leq\mu(C,P)$.  Also, it is not hard to see that
$\mu(C,P)<\infty$ if and only if the characteristic is 0.  (Over
$\mathbb C$, Buchweitz and Greuel, generalizing work of Bassein, name
$\mu(C,P)$ the Milnor number in Def.~1.1.1, p.~244, \cite{BG}, and
prove, in Thm.~4.2.2, p.~258, that, when $C$ degenerates, $\mu(C,P)$
increases by the number of vanishing cycles.)

In characteristic 0,  Buchweitz and Greuel
\cite{BG}, Prop.~1.2.1, p.~246,  prove
 \[
 \mu(C,P)=2\delta(C,P)-r(C,P)+1,
 \]
which extends the Milnor--Jung formula for plane curves.  Now,
$\lambda(C,P)\leq\mu(C,P)$.  So
 $$
 \lambda(C,P)\leq 2\delta(C,P)-r(C,P)+1\eqno(2.1.2)
 $$
 in characteristic 0.  For an upper bound in positive characteristic,
 see Proposition 4.4.
\end{subsct}

\begin{proposition}[Rim]
 Let $C$ be a curve, and $P\in C$ a closed point.  If $C$ is a
 complete intersection at $P$, then $\tau(C,P)=\lambda(C,P)$.
 \end{proposition} \begin{proof} Over $\mathbb C$, the assertion
 follows directly from \cite{BG}, Lem.~1.1.2, p.~245 and Cor.~6.1.6,
 p.~268.  In arbitrary characteristic, the assertion follows directly
 from two formulas buried in the middle of p.~269 in Rim's paper
 \cite{R}.  The first formula says that $\tau$ is equal to the length
 of the torsion submodule of $\Og^1_C$.  A cleaner version of the
 proof, which is based on local duality, was given by Pinkham
 \cite{Pi}, p.~76.  The formula itself was originally proved when $C$
 is irreducible by Zariski, \cite{Z}, Thm.~1, p.~781.  The second
 formula says that this length is equal to $\lambda(C,P)$; here is
 another version of the proof of this formula.

Since the invariants in question are local, we may complete $C$ and then
normalize it off $P$.  Thus we may assume that $C$ is projective and
that $P$ is its only singularity.

The torsion submodule of $\Og^1_C$ is equal to the kernel of the class
map $\gamma\:\Og^1_C \to\w_C$ since $\w_C$ is torsion free.  However,
$\lambda(C,P):=\ell(\Cok(\gamma_P))$.  Hence it suffices to prove
$\chi(\Og^1_C)=\chi(\w_C)$.

Let $\c N$ be the conormal sheaf of $C$ in its ambient projective space,
$X$ say, and set $\c M:=\Og^1_X|C$.  Since $C$ is a local complete
intersection, $\c N$ is locally free and we have an exact sequence of the
form
 \[ 0\to\c  N\to\c  M\to\Og^1_C\to 0.  \]
 So $\chi( \Og^1_C)=\chi(\c M)-\chi(\c N)$.  Hence, by Riemann's theorem,
 \[ \chi(\Og^1_C)=\deg\c  M-\deg\c  N+(\rk\c  M-\rk\c  N)
 \chi(\c O_C)= \deg\c M-\deg\c N+\chi(\c O_C).  \]

 On the other hand, $\w_C=\det(\c M)\ox(\det\c N)^*$ by \cite{Ha},
Thm.~7.11, 
p.~245. So, again   by
 Riemann's theorem,
 \[
 \chi(\w_C)=\deg(\det\c M)-\deg(\det\c N)+\chi(\c O_C).
 \]
Now, $\deg(\det\c M)=\deg\c M$ and $\deg(\det N)=\deg\c N$.  Hence
$\chi(\Og^1_C)=\chi(\w_C)$.
 \end{proof}

\begin{subsct} \emph{Global invariants.}
Let $C$ be a projective curve, $\I n\:\wt C\to C$ the normalization map,
and $\I n^\#\:\c O_C\to\I n_*\c O_{\wt C}$ the associated map.

If $C$ is smooth at a closed point $P$, then the local invariants
$\delta(C,P)$, $\tau(C,P)$, $\lambda(C,P)$, and $\iota(C,P)$ all
vanish.  So it makes sense to set
\begin{alignat*}{2}
 \delta(C)&:=\sum_{P\in C}\delta(C,P),&\quad
 \lambda(C)&:=\sum_{P\in C}\lambda(C,P),\\
 \tau(C)&:=\sum_{P\in C}\tau(C,P),&\quad
 \iota(C)&:=\sum_{P\in C}\iota(C,P).
 \end{alignat*}

Let $r(C)$ denote the number of irreducible components of $C$.

Recall that the  \emph{arithmetic genus} and the \emph{geometric genus}
 are defined by the formulas:
 $$p_a(C):=1-\chi(\c O_C) \text{ \ and \ }p_g(C):=\h^1(\c O_{\wt C}).$$
 Extracting Euler characteristics from the short exact sequence
 \[
 0\to\c O_C\to\I n_*\c O_{\wt C}\to\Cok(\I n^\#)\to 0
 \]
 yields Clebsch's formula
	$$p_g(C)=p_a(C)-\delta(C)+r(C)-1.\eqno(2.3.1)$$
 
Suppose $C$ is connected.  Then $r(C)-1\le\sum_P(r(C,P)-1)$.  In
charactersitic 0, therefore, (2.1.2) yields
	$$\lambda(C)\le 2\delta(C)-r(C)+1.\eqno(2.3.2)$$
\end{subsct}

\begin{proposition}
Let $A$ and $B$ be (reduced) plane curves of degrees $a$ and $b$ with no
common components.  Set $C:=A\cup B$.  Then
 $$ \tau(A)+\tau(B)+ab \le \tau(C) ,$$
 with equality if $A$ and $B$ are transverse.
 \end{proposition}
 \begin{proof}
 If $A$ and $B$ are transverse, then $\tau(C,P)=1$ for $P\in A\cap B$.
There are $ab$ such $P$.  Hence
        $$\tau(A)+\tau(B)+ab=\tau(C).$$

By the theorem of transversality of the general translate for projective
space \cite{K}, Cor.~11, p.~296, there is a dense open subset of
automorphisms $g$ of the plane such that the translate $A^g$ is
transversal to $B$.  Set $C_g:=A^g\cup B$.  Then, by the preceding case,
                         $$\tau(A^g)+\tau(B)+ab=\tau(C_g).$$

The function $g\mapsto\tau(C_g)$ is upper semi-continuous.  Indeed,
$\tau(C_g)=\lambda(C_g)$ by Proposition 2.2.  Furthemore,
$g\mapsto\lambda(C_g)$ is upper semi-continuous, because $\lambda(C_g)$
is the length, on the fiber over $g$, of the restriction of the cokernel
of a map between coherent sheaves on the total space of the $C_g$,
namely, the relative class map.

Hence $\tau(C_g)\le \tau(C)$.  But $\tau(A^g)=\tau(A)$ since $A^g$ and
$A$ are isomorphic.  Therefore, the asserted bound holds.
 \end{proof}

\section{The Hironaka--Noether bound}

\begin{proposition}[Hironaka--Noether bound]
 Let $A$ be a reduced Noetherian local ring of dimension $1$ and
multiplicity $e\ge2$.  Let $B$ be the blowup of $A$ at its maximal ideal
$\mathbf m$.  Then the length of the $A$-module $B/A$ satisfies the
following inequality:
 \[
 \ell(B/A)\le e(e-1)/2.
 \]
Furthermore, equality holds if and only if $A$ has embedding dimension
$2$.
\end{proposition}

\begin{proof} 
 Set $k:=A/\mathbf m$.  Let's first reduce the question to the case
 where $k$ is infinite; we'll use a well-known trick, found for instance
 in \cite{M}, p.~114.  So, let $x$ be an indeterminate, $A[x]$ the
 polynomial ring, and $\mathbf p$ the extension of $\mathbf m$.  Set
 $A(x):=A[x]_{\mathbf p}$.  Then $A(x)$ is a reduced Noetherian local
 ring of dimension~$1$.  Its maximal ideal is the extension $\mathbf m
 A(x)$, and its residue field is the infinite field $k(x)$.

In addition, $A(x)$ is flat over $A$.  Hence, the multiplicity of $A(x)$
is also $e$, and the blowup of $A(x)$ at its maximal ideal is $B\ox_A
A(x)$.  Also,
 $$\ell\bigl((B\ox_A A(x))/A(x)\bigr)=\ell\bigl((B/A)\ox_A A(x)\bigr)
 =\ell(B/A).$$
 Therefore, replacing $A$ by $A(x)$, we may assume $k$ is infinite.

Since $k$ is infinite and $A$ is reduced and of dimension 1, there is an
$f\in \mathbf m$ such that the equation $B=A[\mathbf m/f]$ holds in the
total ring of fractions of $A$.  Note that
	$$\mathbf mB=f(1/f)\mathbf mB\subseteq fB;$$
 whence $\mathbf mB=fB$.  It follows that, for every $i\geq 0$, we have
 $$\ell(\mathbf m^iB/\mathbf m^{i+1}B)=e.\eqno(3.1.1)$$

For each $i\geq 0$, form the  $A$-module  
\[
 V_i:=\mathbf m^iB/(\mathbf m^i+\mathbf m^{i+1}B).
\]
Then $V_i$ is the cokernel of the natural map
\[
\mathbf m^i/\mathbf m^{i+1}\to \mathbf m^iB/\mathbf m^{i+1}B.
\]
 Hence, we get
$$
 \ell(V_i)\geq \ell(\mathbf m^iB/\mathbf m^{i+1}B)
	-\ell(\mathbf m^i/\mathbf m^{i+1}). \eqno(3.1.2)
 $$

Let's now prove that,  for some integer $q\ge 0$, we have
$$
e-1=\ell(V_0)>\ell(V_1)>\cdots>\ell(V_{q-1})
>\ell(V_q)=\ell(V_{q+1})=\cdots=0.
\eqno(3.1.3)
$$
 Indeed, first observe that
\[
 \ell(V_0)=\ell(B/\mathbf mB)-\ell(A/(\mathbf mB\cap A)).
 \]
 Now,  $\ell(B/\mathbf mB)=e$ by (3.1.1).  Also, $\mathbf mB\cap A=\mathbf
m$.  So $\ell(V_0)=e-1$.  

 Next, notice that, for each $i\geq 0$, multiplication by $f$ induces a
map
\[
\begin{CD}
h_i\:V_i @>\times f>> V_{i+1}.
\end{CD}
\]
 This map $h_i$  is surjective because $\mathbf mB=fB$.  Moreover,
$\Ker(h_i)=0$ if and only if \[
  \mathbf m^iB\cap (1/f)(\mathbf m^{i+1}+\mathbf m^{i+2}B)\subseteq
\mathbf m^i+\mathbf m^{i+1}B.
 \]

However, $\mathbf m^{i+1}+\mathbf m^{i+2}B\subseteq \mathbf m^{i+1}B$.
Also $(1/f)\mathbf m^{i+1}B=\mathbf m^iB$ since  $\mathbf mB=fB$. 
Hence, $\Ker(h_i)=0$ if and only if 
\[
(1/f)(\mathbf m^{i+1}+\mathbf m^{i+2}B)
\subseteq \mathbf m^i+\mathbf m^{i+1}B.
\]

Of course, we have
\[
 (1/f)(\mathbf m^{i+1}+\mathbf m^{i+2}B)=(1/f)\mathbf m(\mathbf
m^i+\mathbf m^{i+1}B).
 \]
 Since $B=A[\mathbf m/f]$, it follows that $\Ker(h_i)=0$ if and only if
$\mathbf m^i+\mathbf m^{i+1}B$ is a $B$-module; that is, if and only if
 $$ \mathbf m^i+\mathbf m^{i+1}B=\mathbf m^iB+\mathbf m^{i+1}B
 =\mathbf m^iB. \eqno(3.1.4)
 $$
 Therefore, $h_i\: V_i\to V_{i+1}$ is injective if and only if $V_i=0$.
Since $h_i$ is surjective, $\ell(V_i)\ge\ell(V_{i+1})$; moreover, if
equality holds, then $h_i$ is bijective, and therefore $V_i=0$.  Thus
(3.1.3) holds for some  $q$. 

Next, let's prove that, for all $j\ge 0$, we have
 $$ \mathbf m^q+\mathbf m^{q+j}B=\mathbf m^qB. \eqno(3.1.5)
 $$
 This equation is trivial for $j=0$.  Now, given $j\ge0$, suppose
(3.1.5) holds.  Since (3.1.3) holds, $V_{q+j}=0$; so (3.1.4) holds for
$i:=q+j$.  Hence, we have
 \[ \mathbf m^q+\mathbf m^{q+j+1}B
	=\mathbf m^q+\mathbf m^{q+j}+\mathbf m^{q+j+1}B
	=\mathbf m^q+\mathbf m^{q+j}B=\mathbf m^qB.
 \]
 Thus, by induction, (3.1.5) holds for all $j\ge 0$.

Let's now improve (3.1.5) by showing it implies that
 $$ \mathbf m^q=\mathbf m^qB. \eqno(3.1.6)
 $$
 Indeed, the $A$-module $B/\mathbf m^q$ has finite length.  Hence it is
annihilated by $\mathbf m^{q+j}$ for some $j\ge 0$; in other words,
$\mathbf m^{q+j}B\subseteq \mathbf m^q$.  Thus (3.1.5) yields (3.1.6).

We can now prove the first assertion.  Indeed, owing to (3.1.6), the
sequence
\[
 0\to A/\mathbf m^q\to B/\mathbf m^qB\to B/A\to0
\]
 is exact.  Filter the first term by $\mathbf m^i/\mathbf m^q$ for
$i=0,\dotsc, q$, and the second by $\mathbf m^iB/\mathbf m^qB$.  Then we
get
$$
\ell(B/A)=\sum_{i=0}^{q-1}\bigl(\ell(\mathbf m^iB/\mathbf m^{i+1}B)
 -\ell(\mathbf m^i/\mathbf m^{i+1})\bigr).\eqno(3.1.7)
$$
 Now, (3.1.3) yields $\ell(V_i)\le(e-1-i)$ and $q\le e-1$.  Hence
(3.1.2) yields
 $$
 \ell(B/A)\leq\sum_{i=0}^{q-1}\ell(V_i)
 \leq\sum^{e-2}_{i=0}(e-1-i)=e(e-1)/2.\eqno(3.1.8)
$$
 Thus the first assertion is proved. 

To prove the second assertion, first assume $\ell(B/A)=e(e-1)/2$.  Then
the equalities hold in (3.1.8).  So equality holds in (3.1.2), and
$\ell(V_i)=e-1-i$ for $0\le i\le e-1$.  Hence (3.1.1) yields
$\ell(\mathbf m^i/\mathbf m^{i+1})=i+1$.  In particular, $\ell(\mathbf
m/\mathbf m^2)=2$.

Conversely, assume $\ell(\mathbf m/\mathbf m^2)=2$. Then $\mathbf m$ is
generated by two elements.  So $\mathbf m^i$ is generated by at most
$i+1$ elements for all $i\geq 0$; whence,
  $$\ell(\mathbf m^i/\mathbf m^{i+1}) \leq i+1.\eqno(3.1.9)$$

Together, (3.1.1) and (3.1.6) and (3.1.9) yield
	$$e=\ell(\mathbf m^qB/\mathbf m^{q+1}B)
	=\ell(\mathbf m^q/\mathbf m^{q+1}) \leq q+1.$$
Therefore,  (3.1.7) and (3.1.1)  and (3.1.9) yield
\[
 \ell(B/A)=\sum_{i=0}^{q-1}\bigl(e
	 -\ell(\mathbf m^i/\mathbf m^{i+1})\bigr)
	\ge \sum_{i=0}^{e-2}\bigl(e-1-i)=e(e-1)/2.
\]
Since $\ell(B/A)\le e(e-1)/2$ by (3.1.8), equality holds.
 \end{proof}

\section{Infinitely near points}

\begin{subsct} \emph{Infinitely near points.} Let $X$ be a 
smooth scheme of dimension 2 or more.  An infinite sequence
$P,P',P'',\dotsc$ is said to be a {\it succession of infinitely near
points of} $X$ if $P$ is a closed point of $X$, if $P'$ is a
closed point of the exceptional divisor $E'$ of the blowup $X'$ of
$X$ at $P$, if $P''$ is a closed point of the exceptional divisor
$E''$ of the blowup $X''$ of $X'$ at $P'$, and so forth.

In this case, whenever $m\le n$, then $P^{(n)}$ is said to be {\it
infinitely near to $P^{(m)}$ of order} $n-m$.  In addition, $P^{(n)}$ is
said to be {\it proximate\/} to $P^{(m)}$ if $m<n$ and if $P^{(n)}$ lies
on the proper (or strict) transform of $E^{(m+1)}$ on $X^{(n)}$; given
$n$, denote the number of these $P^{(m)}$ by $i(P,P^{(n)})$.  Note that
$i(P,P^{(n)})=0$ if and only if $n=0$.

Let $C\subset X$ be a curve.  Let $C^{(n)}$ be the proper transform of
$C$ on $X^{(n)}$.  Denote by $e(C,P^{(n)})$, by $\delta(C,P^{(n)})$, and
by $r(C,P^{(n)})$ the multiplicity, the $\delta$-invariant, and the
number of branches of $C^{(n)}$ at $P^{(n)}$; by convention, these
numbers are 0 if $C^{(n)}$ does not contain $P^{(n)}$.  Similarly, given
a branch $\Gamma$ of $C$ at $P$, denote by $e(\Gamma,P^{(n)})$ and by
$\delta(\Gamma,P^{(n)})$ the multiplicity and the $\delta$-invariant at
$P^{(n)}$ of the proper transform of $\Gamma$.

Note that $P^{(n)}$ determines its predecessors $P,P',\dotsc,P^{(n-1)}$,
but not its successors $P^{(n+1)},P^{(n+2)},\dotsc$; the latter vary
with the particular succession through $P^{(n)}$.  Call $P^{(n-1)}$ the
{\it immediate predecessor\/} of $P^{(n)}$.  Denote the set of all
predecessors of $P^{(n)}$, including $P^{(n)}$ and $P$, by
$[P,P^{(n)}]$.  Denote the set of all possible successors $Q$ of
$P^{(n)}$, including $P^{(n)}$, by $N(P^{(n)})$; denote the subset of
those $Q$ proximate to $P^{(n)}$ by $N^*(P^{(n)})$.

\end{subsct}

\begin{lemma}
 Let $X$ be a smooth scheme of dimension $2$ or more, $C\subset X$ a
curve, and $P\in C$ a closed point.  Then
 \[
 \sum_{Q\in N(P)}e(C,Q)\bigl(e(C,Q)-2+i(P,Q)\bigr)
 \geq 2\delta(C,P)-r(C,P),
 \]
 with equality if and only if the embedding dimension of $C$ at $P$ is
   $1$ or $2$.
 \end{lemma}
 \begin{proof}
 The sum in question is well defined.  Indeed, if $Q$ lies off the
proper transform of $C$, then $e(C,Q)=0$.  Of the remaining $Q$, all but
finitely many are such that $e(C,Q)=1$ and $i(P,Q)=1$ by the theorem of
embedded resolution of singularities.

Let $t(C,P)$ be the greatest order of a $Q\in N(P)$ such that either
$e(C,Q)>1$ or $e(C,Q)=1$ and $i(P,Q)>1$.  However, if no such $Q$
exists, set $t(C,P):=-1$.

Suppose $t(C,P)=-1$.  Then, for every $Q\in N(P)\setminus P$, either
$e(C,Q)=0$ or $e(C,Q)=1$ and $i(P,Q)=1$; moreover, $e(C,P)=1$ and
$i(P,P)=0$.  Hence the sum in question is equal to $-1$.  Moreover,
$\delta(C,P)=0$ and $r(C,P)=1$; also the embedding dimension of $C$ at
$P$ is 1.  Hence the assertion holds in this case.

Proceed by induction on $t(C,P)$.  So suppose $t(C,P)\ge0$.  Let $X'$ be
the blowup of $X$ at $P$, and $C'$ the proper transform of $C$.  Say
$P'_1,\dotsc,P'_n\in C'$  lie over $P$.

Fix $j$.  If $t(C',P'_j)=-1$, then $t(C',P'_j)<t(C,P)$.  Now, take $Q\in
N(P'_j)$; say $Q$ is of order $m$.  Then $Q\in N(P)$ with order $m+1$.
Moreover, $e(C',Q)=e(C,Q)$.  Also
	$$i(P'_j,Q)=\begin{cases}
		i(P,Q),& \text{if  $Q$ is not proximate to $P$;}\\
		i(P,Q)-1,& \text{if  $Q$ is  proximate to $P$.}
	\end{cases}$$
 Therefore, if $t(C',P'_j)\ge0$, then again  $t(C',P'_j)<t(C,P)$.

So the induction hypothesis and the above formulas for $e(C',Q)$ and
$i(P'_j,Q)$ yield
 \begin{multline*}
 \sum_{Q\in N(P'_j)}e(C,Q)\bigl(e(C,Q)-2+i(P,Q)\bigr)
 -\sum_{Q\in N(P'_j)\cap N^*(P)}e(C,Q)\\
 \geq2\delta(C',P'_j)-r(C',P'_j),\label[(4.2.1)]
 \end{multline*}
 with equality if the embedding dimension of $C'$ at $P'_j$ is at
most 2.  The latter holds, of course, if the embedding dimension of $C$
at $P$ is at most 2.

Let $\delta$ be the colength of $\c O_{C,P}$ in its blowup.  By
Proposition 3.1,
	$$e(C,P)(e(C,P)-1)\ge2\delta,\eqno(4.2.2)$$
 with equality  if and only if the embedding dimension of $C$
at $P$ is at most 2.  Moreover,
 $$
 \delta(C,P)=\sum_{j=1}^n\delta(C',P'_j)+\delta.\eqno(4.2.3)
 $$

 Sum the inequalities in (4.2.1) over $i$, and use (4.2.2) and (4.2.3).
We get
\begin{align*}
 \sum_{Q\in N(P)}e(C,&Q)\bigl(e(C,Q)-2+i(P,Q)\bigr)\\
 &=e(C,P)(e(C,P)-2)
   +\sum_{j=1}^n\sum_{Q\in N(P'_j)}e(C,Q)\bigl(e(C,Q)-2+i(P,Q)\bigr)\\
 &\ge2\delta -e(C,P)+\sum_{j=1}^n\bigg(\sum_{Q\in N(P'_j)\cap N^*(P)}
  e(C,Q)+2\delta(C',P'_j)-r(C',P'_j)\biggr)\\
 &=2\delta(C,P)-r(C,P)-e(C,P)+\sum_{Q\in N^*(P)}e(C,Q).
\end{align*}
 with equality if and only if the embedding dimension of $C$ at $P$ is
at most 2.  However, the last two terms cancel by the proximity
equality; see \cite{D}, Formula (2.18), p.~27, for example.  Thus the
assertion holds.
 \end{proof}

\begin{lemma}
 Let $X$ be a smooth scheme of dimension $2$ or more in characteristic
$p>0$.  Let $C\subset X$ be a  curve, and $P\in C$ a closed point.
Given a branch $\Gamma$ of $C$ at $P$, let $Q(\Gamma)$ be the point
infinitely 
near to $P$ of least order such that $p\nmid e(\Gamma,Q(\Gamma))$.  Then
	$$\lambda(C,P)\le 2\delta(C,P)-r(C,P)
	+\sum_\Gamma v(\Gamma,P) \text{ where }
	v(\Gamma,P):=\sum_{R\in[P,Q(\Gamma)]}e(\Gamma,R).$$
 \end{lemma}
\begin{proof}
 Let $\I n\:\wt C\to C$ be the normalization map,  $d\,\I n\:
\Og^1_C\to \I n_*\Og^1_{\wt C}$ its differential.  Set
  $$I:=\IM((d\,\I n)_P)\subseteq (\I n_*\Og^1_{\wt C})_P \text{ and }
	\wt I:=(\I n_*\c O_{\wt C})_PI\subseteq (\I n_*\Og^1_{\wt C})_P;$$
 so $I$ is an $\c O_{C,P}$-submodule, and $\wt I$ is the $(\I n_*\c
O_{\wt C})_P$-submodule $I$ generates.  Take an $f\in I$ so that $\wt
I=(\I n_*\c O_{\wt C})_Pf$.  Then $\wt I/(\c O_{C,P}f)\cong (\I n_*\c
O_{\wt C})_P/\c O_{C,P}$.  Hence
	$$\ell(\wt I/I)\leq\delta (C,P).\eqno(4.3.1)$$
 Now, $\I n^*\Og^1_C\to\Og^1_{\wt C}\to\Og^1_{\wt C/C}\to0$ is exact.  So
the Chinese remainder theorem yields
 $$
 (\I n_*\Og^1_{\wt C})_P/\wt I=
 \bigoplus_{\wt P\in\I n^{-1}P}(\Og^1_{\wt C/C})_{\wt P}.\eqno(4.3.2)
 $$

Fix a branch $\Gamma$ of $C$ at $P$, and set $v:=v(\Gamma,P)$.  Say
$\Gamma$ corresponds to $\wt P\in\I n^{-1}P$.  Below, we'll find an
$f\in\c O_{C,P}$ of order $v$ at $\wt P$.  Now, $p\nmid v$.  Hence the
derivative of $f$ with respect to any local parameter of $\wt C$ at $\wt
P$ has order $v-1$.  So $\ell((\Og^1_{\wt C/C})_{\wt P})\le v-1$.

Therefore, Equation (4.3.2) yields
$$
 \ell\bigl((\I n_*\Og^1_{\wt C})_P/\wt I\bigr)
 \leq\sum_{\Gamma}(v(\Gamma,P)-1)=-r(C,P)
 +\sum_{\Gamma}v(\Gamma,P). \eqno(4.3.3)
$$
On the other hand, Equation (2.1.1) yields
 \[
 \lambda(C,P)=\delta(C,P)+\ell((\I n_*\Og^1_{\wt C})_P/I)
  =\delta(C,P)+\ell(\wt I/I)
 +\ell\bigl((\I n_*\Og^1_{\wt C})_P/\wt I\,\bigr).
 \]
Hence, Inequalities (4.3.1) and (4.3.3) yield the assertion, given
the existence of an $f$.

To find an $f$, let $X'$ be the blowup of $X$ at $P$, and $C'$ the
proper transform of $C$.  Say $P'\in C'$ is the image of $\wt P$.  Let
$y_1, \dotsc ,y_m$ be generators of the maximal ideal $\mathbf m_{C,P}$.
Rearranging the $y_i$, we may assume $y_1$ generates the extension
$\mathbf m_{C,P}\c O_{C',P'}$.  Then the order of $y_1$ at $\wt P$ is
$e(\Gamma,P)$.  So, if $p\nmid e(\Gamma,P)$, that is, if $Q=P$, take
$f:=y_1$.

Proceed by induction on the order $n$ of $Q/P$.  Suppose $n>0$.  Then
the order of $Q/P'$ is $n-1$.  Say $y_i=z_iy_1$ where $z_i\in\c
O_{C',P'}$.  Let $a_i$ be the value $z_i$ takes at $P'$.  Then
$y_1,z_2-a_2, \dotsc ,z_m-a_m$ are generators of the maximal ideal
$\mathbf m_{C',P'}$.

By induction, we may assume that a certain scalar linear combination
	$$f':=b_1y_1 + b_2(z_2-a_2) +\dotsb+ b_m(z_m-a_m)$$
 has order $v(\Gamma,P')$ at $\wt P$.  Then $f'y_1$ has order
$v(\Gamma,P)$ at $\wt P$.  Furthermore, $f'y_1$ is a scalar linear
combination of the $y_i$.  So take $f:=f'y_1$.
 \end{proof}

\begin{proposition} 
 Let $X$ be a smooth scheme of dimension $2$ or more in characteristic
 $p\ge0$.  Let $C\subset X$ be a  curve, and $P\in C$ a closed
 point.  If $p=0$, then
 \[
 \lambda(C,P)\le 1+ \sum_{Q\in N(P)}e(C,Q)\bigl(e(C,Q)-2+i(P,Q)\bigr).
 \]

Suppose $p>0$.  For each $Q\in N(P)$, set $\epsilon(C,Q):=0$ if $Q\neq
P$ and if $e(C,R)\le1$ where $R$ is the immediate predecessor of $Q$;
otherwise, set $\epsilon(C,Q):=1$.  Then
 \[
 \lambda(C,P)
 \le \sum_{Q\in N(P)}e(C,Q)\bigl(e(C,Q)-2+i(P,Q)+\epsilon(C,Q)\bigr),
 \]
 \end{proposition}

 \begin{proof}
 If $p=0$, then the asserted bound follows directly from (2.1.2) and
Lemma 4.2.

Suppose $p>0$.  Fix $Q\in N(P)$.  Notice, as $\Gamma$ ranges over all
the branches of $C$ at $P$,
	$$\sum_{\Gamma}e(\Gamma,Q)=e(C,Q). \eqno(4.4.1)$$

Fix a $\Gamma$, and suppose $Q$ is the point of least order such that
$p\nmid e(\Gamma,Q)$.  Let $R\in[P,Q]$.  If $R\neq Q$, then $p\mid
e(\Gamma,R)$, and so $e(\Gamma,R)>1$.  Hence $\epsilon(C,R):=1$ for all
$R\in[P,Q]$.

It now follows from Lemma 4.3 and Formula (4.4.1) that
	$$\lambda(C,P)\le 2\delta(C,P)-r(C,P)
	+\sum_{Q\in N(P)}e(C,Q)\epsilon(C,Q).$$
 Hence Lemma 4.2 yields the asserted bound. 
 \end{proof}

\section{Foliations}

\begin{subsct} \emph{Foliations.} Let 
$X$ be a scheme, $\c L$ an invertible sheaf, and $\eta\:\Og^1_X\to\c L$
a nonzero map.  Then $\eta$ will be called a \emph{(singular
  one-dimensional)  foliation} of $X$.

Let $S\subseteq X$ be the zero scheme of $\eta$, that is, the closed
subscheme whose ideal $\c I_{S/X}$ is the image of the induced map
$\Og^1_X\ox\c L^{-1}\to\c O_X$.  Then $S$ will be called the
\emph{singular locus} of $\eta$.

Let $C\subseteq X$ be a closed curve.  Suppose for a moment (1) that
$C\cap S$ is finite and (2) that the restricted map $\eta|C$ factors
through the standard map $\sigma\:\Og^1_X|C\to\Og^1_C$, in other words,
that there is a commutative diagram
 $$
 \begin{CD}
 \Og^1_X @>\eta>> \c L\\
 @V VV @VVV\\
 \Og^1_C @>\mu>> \c L|C
 \end{CD}
 \eqno(5.1.1)
 $$
  Then $C$ will be called a \emph{leaf} of $\eta$.

Notice the following.  Assume $X$ is smooth.  Let $P\in X-S$ be a
closed point, and $\eta^*\:\c L^*\to \T_X$ the dual map.  Then the
image of $\eta^*(P)$ is a one-dimensional vector subspace, $F(P)$ say,
of the fiber $\T_X(P)$.  Moreover, if $C$ is a leaf and if $P$ is a
simple point of $C$, then $F(P)\subseteq \T_C(P)$.

Conversely, assume $C\cap S$ is finite, and let $U\subseteq C-S$ be a
dense open subset.  Let's prove that, if $F(P)\subseteq \T_C(P)$ for
every simple point $P\in U$, then $C$ is a leaf.

Indeed, let $\c K$ be the kernel of $\sigma\:\Og^1_X|C\to\Og^1_C$, and
$\kappa\:\c K\to\c L|C$ the restriction of $\eta|C$ to $\c K$.  It
follows from the hypothesis that $\kappa(P)=0$ for every simple point
$P\in U$.
So, since $U$ is dense in $C$, the image of $\kappa$ has finite
support.  Now, $C$ is reduced and $\c L|C$ is invertible.  Hence
$\kappa=0$. So there is a map $\mu\:\Og^1_C\to\c L|C$ making the
diagram (5.1.1) commute.  Thus $C$ is a leaf.
 \end{subsct}

\begin{proposition}
 Let $X$ be a scheme, $C\subseteq X$ a projective curve,
 $\eta\:\Og^1_X\to\c L$ a foliation, and $S$ its singular locus. If $C$
 is a leaf of $\eta$, then
  \begin{align*}
 2p_a(C)-2-\deg(\c L|C)&=\lambda(C)-\deg(C\cap S)\\
	&\le\lambda(C)-\iota(C).
  \end{align*}
 \end{proposition}
 \begin{proof}
 Form the standard exact sequence
	$$
	0\to\c I_{(C\cap S)/C}\to\c O_C \to\c O_{C\cap S}\to 0.
 $$
 Twist it by $\c L$, and take Euler characteristics; we get
  $$\chi(\c I_{(C\cap S)/C}\ox\c L) = \chi(\c L|C)-\chi(\c L|(C\cap S)).$$
 Use Riemann's theorem to evaluate $\chi(\c L|C)$.  Then we get
  $$\chi(\c I_{(C\cap S)/C}\ox\c L)
	= \deg(\c L|C)+1-p_a(C)-\chi(\c L|(C\cap S)).\eqno(5.2.1)$$
  
Since $C$ is a leaf, there is a map $\mu\:\Og^1_C\to\c L|C$ making the
diagram (5.1.1) commute.  Since $S$ is the singular locus of $\eta$, the
image $\IM(\eta)$ is equal to $\c I_{S/X}\ox\c L$.  Hence
 $$\IM(\mu)=\c I_{(C\cap S)/C}\ox\c L.\eqno(5.2.2)$$

So $\Cok(\mu)=\c L|(C\cap S)$.  However, $\c L$  is invertible.  Hence
      $$\iota(C)\le\chi(\c L|(C\cap S))=\deg(C\cap S).\eqno(5.2.3)$$

On the other hand, $C$ is reduced.  So $\c L|C$ is torsion free.  Hence
$\IM(\mu)$ is equal to $\Og^1_C/\text{torsion}$ because $C\cap S$ is
finite.  In addition, the canonical sheaf $\w_C$ is torsion free.  Hence
the image of the class map $\gamma\:\Og^1_C\to\w_C$ is also equal to
$\Og^1_C/\text{torsion}$.  So
	$$\IM(\gamma)=\IM(\mu).\eqno(5.2.4)$$

Since $\lambda(C)=\chi(\Cok(\gamma))$, it follows that
  $$\lambda(C)=\chi(\w_C)-\chi(\IM(\gamma)).$$
 Now, $\chi(\w_C)=p_a(C)-1$.  
 Hence (5.2.1)--(5.2.4) yield the assertion.
 \end{proof}

\begin{theorem}
Let $X$ be a smooth scheme of dimension $2$ or more in characteristic
$p\ge0$, and $C\subset X$ a projective curve.  Let $P$ range over all the
closed points of $X$.  For each $Q\in N(P)$, set $\epsilon(C,Q):=0$
either (i) if $e(C,Q)=0$, or (ii) if $p=0$, or (iii) if $p>0$, if
$Q\neq P$, and if $e(C,R)=1$ where $R$ is the immediate predecessor of
$Q$; otherwise, set $\epsilon(C,Q):=1$.  Next, set
	$$\ell(C,Q):=e(C,Q)-2+i(P,Q)+\epsilon(C,Q).$$

 $(1)$ Let $\eta\:\Og^1_X\to\c L$ be a foliation, and assume $C$ is a
leaf.  Then
	$$2p_a(C)-2-\deg(\c L|C)
	\le\sum_{P\in X}\sum_{\,Q\in N(P)} e(C,Q)\ell(C,Q).$$

 $(2)$ Let $A\subset X$ be a divisor.  For each $P$ and $Q\in N(P)$, let
$e(A,Q)$ be the multiplicity at $Q$ of the proper transform of $A$ on
the successive blowup of $X$ determined by $Q$.  Assume that
$e(A,Q)\ge\ell(C,Q)$ and that $C$ is a leaf of $\eta\:\Og^1_X\to\c L$.
Then
	$$2p_a(C)-2-\deg(\c L|C)\le (A\cdot C).$$ 
 \end{theorem}

\begin{proof}
 To prove (1), recall that, if $p=0$ and $P$ is a singular point of $C$,
then
 $\iota(C,P)\ge1$.  Hence Theorem 5.2 and Proposition 4.4 yield (1).

To prove (2), note that $(A\cdot C)=\sum_Q e(A,Q)e(C,Q)$ by Noether's
formula; see \cite{D}, Formula (2.17), p.~27, for example.  Hence (1)
yields (2).
 \end{proof}

\section{Projective space}

\begin{theorem}
Let $X:=\IP^n$ with $n\ge2$, and let $C\subset X$ be a closed curve of
degree $d$.
Assume  $d$ is not a multiple of the characteristic.
  Let $\eta\:\Og^1_X\to\c O_X(m-1)$ be a foliation, $S$
its singular locus.  Assume $C$ is a leaf.  Then
	$$2p_a(C)-(d-1)(m-1)\le\lambda(C),$$
  with equality only if $C\cap S$ has degree $m+1$ and lies on a line
  $M$ and either $M\subseteq S$ or $M$ is a leaf.
 \end{theorem}

\begin{proof}
It is well known, and reproved below, that $\deg(C\cap S)$ is at least
the Castelnuovo--Mumford regularity $\reg(C\cap S)$.  In turn,
$\reg(C\cap S)\ge m+1$ owing to \cite{EK}, Cor.~4.5.  So Proposition 5.2
yields the asserted inequality.

Suppose equality holds in the assertion.  Then the above reasoning yields
 $$\deg(C\cap S)=\reg(C\cap S)=m+1.\eqno (6.1.1)$$ It follows, as is
well known and reproved below, that the scheme $C\cap S$ lies on
 a line $M$.

Suppose that $M\not\subseteq S$ and that $M$ is not a leaf.   Then there
is a point  $P$ in $M\setminus S$ at which the tangent ``direction''
$F(P)\subset T_{X,P}$ associated to $\eta$ differs from that $T_{M,P}$ 
associated to $M$; see the end of Subsection 5.1.  Take a hyperplane $H$ 
containing $M$ such that $T_{H,P}\not\supset F(P)$.

Let $\beta\:\Og^1_X|H\to\Og^1_H$ be the natural map, and set
$\xi:=(\beta,\,\eta|H)$, so that
 $$\xi\:\Og^1_X|H\to\Og^1_H\oplus\c O_H(m-1).$$ 
 Set $\zeta:=(\wedge^n\xi)(n+1)$.  Now, $\eta|H$ factors through the
 twisted ideal $\c I_{(H\cap S)/H}(m-1)$.  So $\zeta$ factors through
   $\c I_{(H\cap S)/H}(m)$.
However, $\zeta(P)\ne 0$ because $T_{H,P}\not\supset F(P)$.

Form the zero scheme $Z$ of $\zeta$.  It follows that $\c O_H(Z)=\c
O_H(m)$; 
also, $Z\supset H\cap S$, but $Z\not\ni P$, whence $Z\not\supset M$.  
So $M\cap Z$ is
finite, has degree $m$, and contains $M\cap S$.  But $\deg(M\cap S)\ge
m+1$ because $(M\cap S)\supseteq (C\cap S)$ and because of (6.1.1).  A
contradiction has been reached.  So the proof is now complete, given the
two well-known results.

Let's now derive these two results from Mumford's original work
\cite{Mm}.  Let $W\subset X$ be a finite subscheme.  Take a hyperplane
$H$ that misses $W$.  Then the ideal $\c I_{(H\cap W)/H}$ is trivial, so
it is 0-regular.  Hence, by the last display on p.~102 in \cite{Mm}, the
ideal $\c I_{W/X}$ is $r$-regular with $r:=h^1(\c I_{W/X}(-1))$.  But
$r=\deg W$ owing to the sequence
	$$0\to\c I_{W/X}(s)\to\c O_X(s)\to\c O_W(s)\to0$$
 with $s:=-1$.  Thus $\reg W\le\deg W$.

Suppose $\reg W=\deg W$.  So $\h^1(\c I_{W/X}(\deg W-2))\ne 0$. As 
$h^1(\c I_{W/X}(-1))=\deg W$, it follows that 
$\h^1(\c I_{W/X}(1))=\deg W-2$, by
Display $(\#')$ on p.~102 in \cite{Mm}.  Hence $\h^0(\c I_{W/X}(1))=n-1$
owing 
to the above sequence with $s:=1$.  So $W$ lies on $n-1$ linearly 
independent hyperplanes of $X$, whence on their line of intersection.
 \end{proof}

\begin{corollary} Let $X:=\IP^n$ with $n\ge2$, and $C\subset X$ be a
closed curve of degree $d$.
Assume   $C$ is connected and the characteristic is $0$.
  Let $\eta\:\Og^1_X\to\c O_X(m-1)$ be a foliation.  Assume $C$ is a
  leaf.  Then
 \[
  p_g(C)\leq (m-1)(d-1)/2+(r(C)-1)/2.
 \]
 \end{corollary}
 \begin{proof}
 The assertion results  from  Theorem 6.1, 
 Formula (2.3.1), and  Bound (2.3.2).
\end{proof}

\begin{proposition}
Let $X:=\IP^n$ with $n\ge2$, and $C\subset X$ a closed curve of degree
$d$.
 Let $\eta\:\Og^1_X\to\c O_X(m-1)$ be a foliation,
$S$ its singular locus.  Assume $S$ is finite and $C$ is a leaf.  Then
	$$\lambda(C)\le2p_a(C)-(d-1)(m-1)+m^2+\dotsb+m^n.$$
 \end{proposition}

\begin{proof}
 Since $S$  is finite, it represents the top Chern class of
$(\Og^1_X)^*(m-1)$.  Hence
	$$\deg(S)=1+m+m^2+\dotsb+m^n.$$
 Since  $\deg(S)\ge\deg(C\cap S)$, Proposition 5.2 now yields the
assertion.
 \end{proof}

\begin{corollary}[du Plessis and Wall]
Let $C$ be a (reduced) plane curve of degree $d$.
Assume  $d$ is not a multiple of the characteristic.
 Let $m$ be the least degree of a nonzero polynomial vector field
 $\phi$ annihilating the polynomial defining $C$.  Then $m\le d-1$ and
	$$(d-1)(d-m-1)\le\tau(C).$$
 If the foliation defined by $\phi$ has only finitely many singularities
on $C$, then also $$\tau(C)\le (d-1)(d-m-1)+m^2.$$
 \end{corollary}

\begin{proof}
 Pick homogeneous coordinates $x,y,z$ for the plane $X$.  Say
	$$\phi=\fld fgh\text{ \ and \ } C:u=0$$
 where $f,g,h$ are polynomials in $x,y,z$ of degree $m$ and where $u$ is
 one of degree $d$.  By hypothesis, $\phi u=0$. Also, $\phi\ne 0$; that
is, 
$(f,g,h)\ne 0$.

In any case, $u$ is annihilated by the three Hamilton fields
	$$\Ham uyz,\ \ \Ham uzx,\ \ \Ham uxy.$$
Since  $d$ is not a multiple of the characteristic,
at least two of the three are nonzero.  Hence
$m\le d-1$.

Consider the Euler exact sequence,
 \[
 \begin{CD}
 0@>>>\Og^1_X @>>> \c O_X(-1)^3 @>(x,y,z)>>\c O_X@>>> 0.
 \end{CD}
 \]
The triple $(f,g,h)$ defines a map $\c O_X(-1)^3\to\c O_X(m-1)$.  Let
$\eta$ be its restriction to $\Og^1_X$.

Owing to the exactness, $\eta=0$ if and only if $(f,g,h)=p(x,y,z)$ for
some polynomial $p$.  But, if $p$ exists, then $\phi u=0$ yields
$pdu=0$; whence, $p=0$ because $d$ is not a multiple of the
characteristic.  Since $\phi\ne 0$, necessarily $\eta\ne 0$.  Thus
$\eta$ is a foliation. 

Diagram (5.1.1) exists as $\phi u=0$.  So, if $C\cap S$ is finite, then
$C$ is a leaf.

Since $C$ is plane, $\tau(C)=\lambda(C)$ by Proposition 2.2; also,
$2p_a(C)-2=d(d-3)$ by adjunction.  Therefore, if $S$ is finite, and so
$C$ is a leaf, then the asserted bounds follow from Theorem 6.1 and
Proposition 6.3.

So assume $S$ is infinite.  Let $B\subseteq S$ be the effective divisor
of largest degree, $b$ say.  Then $\c I_{S/X}\subseteq\c O_X(-B)$.  So
$\eta$ factors through a foliation $\eta'\:\Og^1_X\to\c O_X(m-1-b)$,
whose singular locus has $\c I_{S/X}(B)$ as its ideal.  Hence the
singular locus of $\eta'$ is finite.

Set $\c L:=\c O_X(m-1-b)$.  The Euler sequence gives rise to the sequence
 $$\Hom(\c O_X(-1)^3,\c L)\to\Hom(\Og^1_X,\c L)\to\Ext^1(\c O_X,\c L).$$
The third term is equal to $H^1(\c L)$, so vanishes.  Hence $\eta'$
lifts to a polynomial vector field $\phi'$ of degree $m-b$.  Then
$\phi'\ne0$ simply because  $\phi'$ is a lift.

Say $B:w=0$ where $w$ is a polynomial of degree $b$.  Then
$\phi-w\phi'=p\epsilon$ where $p$ is a suitable polynomial and 
	$$\epsilon:=\fld xyz$$
 is the Euler, or radial, vector field.    Now, $\phi u=0$; hence,
  $$-w\phi'u=pdu.\eqno(6.4.1)$$

Let $T$ be a component of $B$.  Say $T:t=0$ where $t$ is a polynomial of
degree $e$.   Suppose $T$ is not a component of $C$.  Then $t\mid w$, but
$t\nmid u$.  So (6.4.1) implies $t\mid p$.  Set $q:=p/t$ and $r:=w/t$.
Then $r\phi'u=-qdu$.  Set $\phi'':=r\phi'+q\epsilon$.  Then $\phi''u=0$.
Moreover, $\phi''\neq 0$ because $\eta\neq 0$.  So $\phi''$ is a nonzero
polynomial vector field of degree $m-e$ annihilating $u$.  But $m-e<m$,
yet $m$ is minimal---a contradiction!  Thus $T$ is a component of $C$.

Suppose $T$ appears in $B$ with multiplicity 2 or more.  Set $r:=w/t^2$.
Since $u$ is reduced, (6.4.1) implies $t\mid p$.  Set $q:=p/t$.  Then
$rt\phi'u=-qdu$.  Set $\phi'':=rt\phi'+q\epsilon$.  Then $\phi''$ is a
nonzero polynomial vector field of degree $m-e$ annihilating $u$.  But
$m-e<m$, yet $m$ is minimal---a contradiction!  Thus $B$ is reduced.

 Set $A:=C-B$ and $a:=d-b$.  Then $A$ is a reduced effective divisor, so
a curve of degree $a$.  And  $a>0$ as   $b\le m \le d-1$.
 Moreover, $A$ is a leaf of $\eta'$, which has finite singular locus.
 As observed above, Proposition 2.2, adjunction, and Theorem 6.1 yield
	$$(a-1)(a-(m-b)-1)\le\tau(A).$$
 Now, $\tau(A)+\tau(B)+ab\le\tau(C)$ by Proposition 2.4.  But
$0\le\tau(B)$.
 Hence
	$$(a-1)(a-(m-b)-1)+ab\le\tau(C).$$
 Now, $m\ge b$; so  $b(a-(m-b)-1)<ab$.  Hence
	$$(a+b-1)(a+b-m-1)<(a-1)(a-(m-b)-1)+ab\le\tau(C).$$
  Since $a+b=d$, the first assertion therefore holds.

As to the second assertion, suppose $\eta$ has only finitely many
singularities on $C$.  But $B\subset C$.  Hence $B=\emptyset$.  So $S$
is finite.  Therefore, as was observed above, the upper bound holds.
\end{proof}

\smallskip
 \centerline{\vrule width 3 true cm height 0.4pt depth 0pt}

\medskip

\small
  \centerline{IMPA,
 Estrada D. Castorina {\sl110, 22460--320} Rio de Janeiro RJ, BRASIL}
  \smallskip

  \centerline{E-mail: {\texttt{esteves@impa.br}}}

 \medskip

  \centerline{Math Dept, Room {\sl 2-278} MIT,
   {\sl77} Mass Ave, Cambridge, MA {\sl02139-4307}, USA}
 \smallskip

  \centerline{E-mail: {\tt kleiman@math.mit.edu}}

\end{document}